\documentclass{ws-adsaa}

\usepackage{xcolor}
\usepackage[verbose]{hyperref}
\hypersetup{colorlinks=false,allbordercolors=blue,pdfborderstyle={/S/U/W 1}}
\def\R{{\mathbb{R}}}

\def\Z{{\mathbb{Z}}}
\def\N{{\mathcal{N}}}

\def\F{{\mathcal{F}}}

\def\W{{\mathcal{W}}}

\DeclareMathOperator\supp{supp} % support

\begin{document}

\markboth{Authors' Names}{Instructions for Typing Manuscripts (Paper's Title)}

%%%%%%%%%%%%%%%%%%%%% Publisher's Area please ignore %%%%%%%%%%%%%%%
%
\catchline{}{}{}{}{}
%
%%%%%%%%%%%%%%%%%%%%%%%%%%%%%%%%%%%%%%%%%%%%%%%%%%%%%%%%%%%%%%%%%%%%

\title{Empirical wavelet frames}

\author{J\'er\^ome Gilles}

\address{Department of Mathematics \& Statistics, San Diego State University, 5500 Campanile Dr\\ San Diego, California, 92182, USA,\\
\email{jgilles@sdsu.edu} \http{jegilles.sdsu.edu}}

\author{Richard Castro}

\address{Department of Mathematics \& Statistics, San Diego State University, 5500 Campanile Dr\\ San Diego, California, 92182, USA,\\
\email{rcastro0091@gmail.com}}

\maketitle

\begin{history}
\received{Day Month Year}
\revised{Day Month Year}
\accepted{Day Month Year}
\published{Day Month Year}
\end{history}

\begin{abstract}
  Due to their adaptive nature, empirical wavelets had several successes in many fields from engineering, science, medical signal/image processing. Recently, a general theoretical framework has been developed in the one-dimensional case, showing the possibility to build empirical wavelets from any classic mother wavelets. Given extensive literature both in theory and applications of classic wavelet frames, it is legitimate to ask about the feasibility of building empirical wavelet frames. We address this question in this paper. We prove several results which provide conditions on the existence of empirical wavelet frames taking into account the above mentioned adaptability.
\end{abstract}


\section{Introduction}\label{sec1}
The empirical wavelet transform in one dimension has been introduced in \cite{Gilles2013} and extended into two dimension in \cite{Gilles2013a}, \cite{ewwt}, \cite{evw}. The concept of empirical wavelets (EW) is to propose a wavelet type transform whose corresponding filter bank is adaptive, i.e. data-driven. The main idea, inspired by the empirical mode decomposition (EMD) \cite{Huang1998}, consists in considering the input signal $f$ as a sum of harmonic modes (i.e. amplitude modulated - frequency modulated components) $f_n$, plus a residue $r$, i.e.
$$f(t)=r(t)+\sum_{n=1}^N f_n(t) \qquad \text{where}\quad f_n(t)=A_n(t)
cos(\omega_n(t)),$$
and $N$ is the number of modes which is either given or has to be estimated. It is also assumed that $A_n(t)\geq 0$ and $\frac{d\omega_n(t)}{dt}\geq 0$ for all $t$ and $n$ to guarantee that the modes $f_n$ are distinct ``enough'' in the Fourier domain. Such decomposition can be easily achieved by assuming that such harmonic modes have compact supports (or are very fast decaying) in the Fourier domain. Thus separating them corresponds to detecting the position of their respective supports in the Fourier domain. Then, wavelet filters are built with respect to each support, providing the expected filter bank used to perform the final decomposition.\\
This transform has proven to be useful for the analysis of non-stationary signals in many fields of science and engineering \cite{Huang2018}, \cite{Huang2019}, \cite{lfm}, \cite{wind}, \cite{ecg}, \cite{speech}, \cite{singlephase}, \cite{seismic}, \cite{EEGEWT}, \cite{klaar2023optimized}, \cite{yan2023wavelet}, \cite{liu2023adaptive}, \cite{mohammadi2023using}, \cite{mohapatra2023gastrointestinal}, \cite{mo2023conditional}, \cite{el2022rescovidtcnnet}, \cite{peng2022effective} just to mention a few.\\
The original EW were designed using a specific Littlewood-Paley type wavelet filter, a general mathematical framework was developed in \cite{cewt}, allowing the use of all classic mother wavelets defined either in the Fourier or the time domain. A general condition is also given in \cite{cewt} to guarantee the existence of the reconstruction formula.\\

Harmonic decompositions have been widely studied through the formalism of frames. We recall that a countable family of functions $\{\psi_n\}$ is a frame of $L^2(\R)$ if for all function $f\in L^2(\R)$, there exists two constants $0<A\leq B<\infty$, such that
$$A\|f\|_{L^2}^2\leq \sum_n |\langle f,\psi_n\rangle|^2\leq B\|f\|_{L^2}^2.$$
If $A=B$, the frame is said tight, and if $A=B=1$ then $\{\psi_n\}$ is a Parseval frame. If only the right inequality holds then the family of function is called a Bessel sequence. Frames have been extensively studied in the literature \cite{Janssen1996}, \cite{Azarmi2017}, \cite{Balazs2011}, \cite{Bhat2018}, \cite{Bownik2012}, \cite{Cabrelli2014}, \cite{Christensen2010}, \cite{Chui1993}, \cite{Doerfler2015}, \cite{Nambudiria2018}, \cite{Sun2002}, especially frames built from wave packet systems, i.e. made from the translation, modulation or dilation (or any combination of these operators) of one or several ``mother'' functions \cite{Labate2002}, \cite{Labate2004}. For instance, classic wavelets correspond to the combination of translations and dilations, while modulations and translations provide Gabor or Weyl-Heisenberg frames.\\

Given that in \cite{cewt}, it is shown that empirical wavelets can be written as special intertwined combination of a modulation, translation and dilation, it is natural to ask the question if such empirical wavelet systems can form frames? We answer this question in the work presented in this paper. The remainder of the paper is organized as follow. In Section~\ref{sec1b}, we set up some notations that will be used throughout the paper. In Section~\ref{sec2}, we recall the formalism of empirical wavelet systems, and define extra notations that will simplify our work. In Section~\ref{sec3}, we briefly recall some results about wave packet systems which will be useful to prove our results. Our main results will be stated and proven in Section~\ref{sec4} where we explore the construction of empirical wavelet frames. Finally, we will conclude and give some perspectives in Section~\ref{sec5}.

\section{Notations}\label{sec1b}
In this section, we recall some classic definitions and notations we will use throughout the paper. All integrals will be considered in the sense of Lebesgue. We will consider the space $L^2(\R)$ equipped with its standard inner product defined by,
$$\forall f,g\in L^2(\R),\qquad \langle f,g\rangle=\int_\R f(x)\overline{g(x)}dx,$$
which induces the standard $L^2$ norm, i.e. $\forall f\in L^2(\R),\|f\|_{L^2}^2=\int_\R |f(x)|^2dx$. The Fourier transform will be defined by, $\forall f\in L^1(\R)\cap L^2(\R)$,
$$\hat{f}(\xi)=(\F f)(\xi)=\int_\R f(x) e^{-2\imath\pi \xi x} dx,$$
and its inverse by 
$$f(x)=(\F^{-1} \hat{f})(x)=(\hat{f})^\vee=\int_\R \hat{f}(\xi) e^{2\imath\pi \xi x} d\xi.$$
The Fourier transform is classically extended to $L^2(\R)$ as a unitary operator.\\

Let us denote, for any arbitrary function $f\in L^2(\R^n)$, the following operators: the translations $T_y:(T_yf)(x)=f(x-y)$ where $y\in\R^n$; the dilations $D_a:(D_af)(x)=|\det a|^{-1/2}f(a^{-1}x)$ where $a\in GL_n(\R)$; and the modulations $E_\nu:(E_\nu f)(x)=e^{2\imath\pi\nu\cdot x}f(x)$ where $\nu\in\R^n$. It is straightforward to check that these operators fulfill the following properties:
\begin{align*}
    \F(E_\nu f) &= T_\nu \widehat{f} & &; &   \F^{-1}(T_y \widehat{f})&= E_y f,\\
     \F(T_y f) &= E_{-y} \widehat{f} & &; & \F^{-1}(E_\nu \widehat{f}) &= T_{-\nu} f,\\
     \F(D_a f) &= D_{a^{-1}} \widehat{f}  & &; & \F^{-1}(D_a \widehat{f}) &= D_{a^{-1}} f.
\end{align*}

\section{Empirical wavelets systems}\label{sec2}
The particularity of empirical wavelets (EW), is in the fact that the filters of the corresponding filter bank are not designed based on a prescribed partitioning of the Fourier domain. Instead, the ``optimal'' partitioning is detected from the spectra of the function we aim to decompose (from a practical point of view, the EW transform is data-driven). In \cite{cewt}, a general formalism is provided for the construction of EW in the one-dimensional case, we recall this formalism and make a few minor modifications to make it more convenient. We first recall how Fourier partitions are defined. Secondly, we provide the formalism to construct EW.

%=================================================================
\subsection{Partitioning of the Fourier Domain}
%=================================================================
In this section, we describe a partitioning of the Fourier line (we only consider the one-dimensional case, $\R$), as in \cite{cewt}, which will be used to construct empirical wavelets.
First, let $n_m, n_M \in \Z$, with $n_m < n_M$, we define two indexing sets as follows 
\begin{align*}
\N &= \{ n_m, \ldots ,0, \ldots , n_M \}, \\
\N^* &= \N \setminus \{ 0 \}.   
\end{align*}

To make them distinguishable, we will denote sets indexed by $\N ^*$ with *, i.e. $S = \{ s_n \}_{n \in \N}$ and $S^* = \{ s_n \}_{n \in \N^*}$.

Let us denote $\mathcal{V} = \{ \nu_n \}_{n \in \N} \subset \overline{\R}$, where the points $\nu_n$ are arbitrary boundary points, and $\forall n\in \Z,\nu_n<\nu_{n+1}$. We adopt the convention $\nu_0 = 0$, which implies that if $n < 0$ then $\nu_n < 0$ and if $n > 0$ then $\nu_n > 0$. In some cases $\nu_0 = 0$ can be excluded, and in such a case we will denote the partition $\mathcal{V}^*$. Such set of boundary points allows us to define Fourier partitions.\\

\begin{figure}[t!]
    \includegraphics[width = \linewidth]{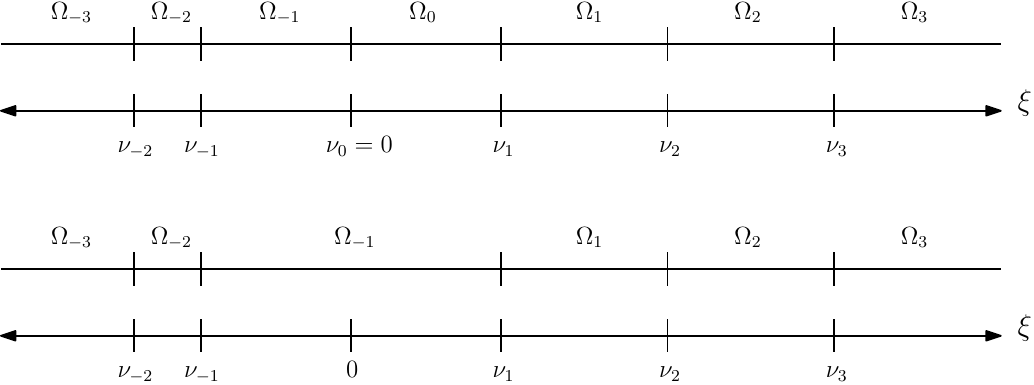}
    \caption{Example of partitions $\mathcal{V}$ (top) and $\mathcal{V}^*$ (bottom) of the Fourier line.}
    \label{fig:4.1}
\end{figure}

\begin{definition} \label{Defn 4.1}
A partition of the Fourier line based on the set of boundary points $\mathcal{V}$ can be defined as the set of intervals, or \textbf{Fourier supports}, of the form 
\begin{align}
\forall n \in \N,\quad \Omega_n = [\nu_n, \nu_{n+1}].         \label{4.3}
\end{align}
\end{definition}

\begin{figure}[t!]
    \includegraphics[width = \linewidth]{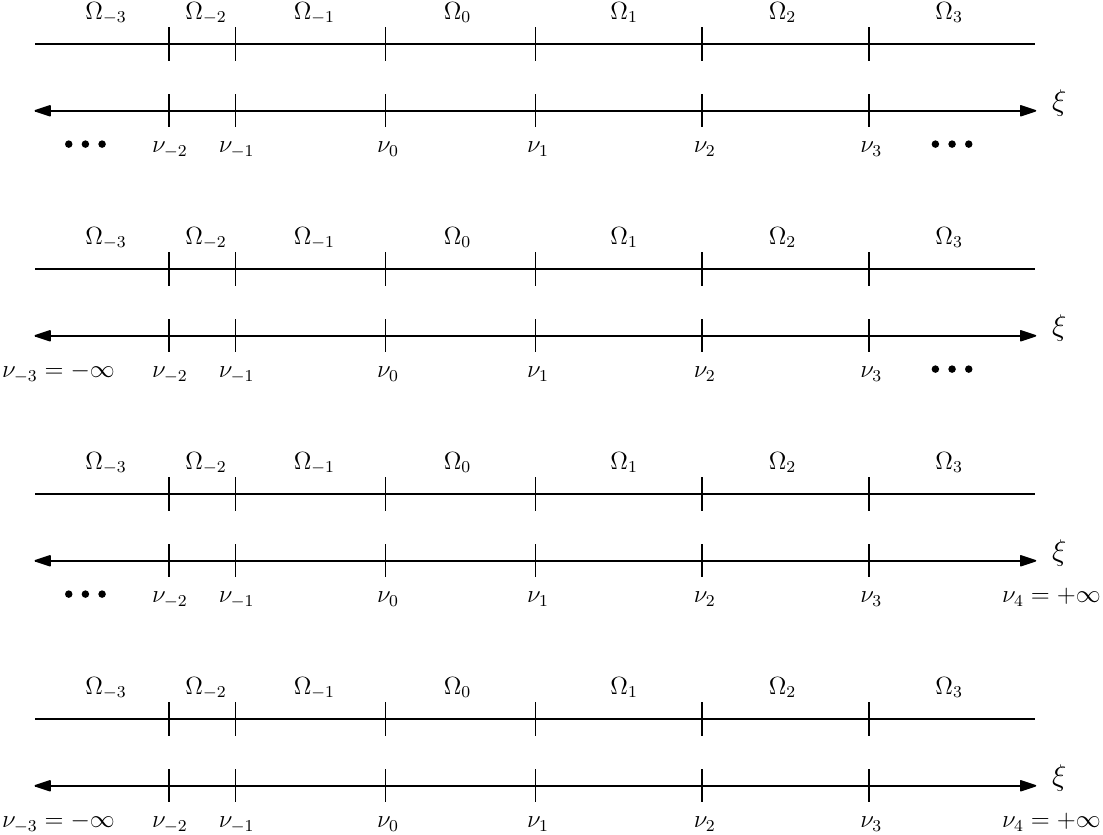}
    \caption{Example of the four types of partitions. From top to bottom: $\Omega^{I}$, $\Omega^{I}_{Lray}$, $\Omega^{I}_{Rray}$, and $\Omega^{F}_{rays}$.}
    \label{fig:4.2}
\end{figure}

In the case of $\mathcal{V}^{*}$, we use the same definition, Definition~\ref{Defn 4.1}, except that $n \in \N^{*}$ and we set $\Omega_{-1} = [\nu_{-1}, \nu_{1}]$. In the remainder of this paper, we will denote either $\Omega = \{ \Omega_n \}_{n \in \N}$ or $\Omega^{*} = \{ \Omega_n \}_{n \in \mathcal{N^*}}$ a given partition of the Fourier line. An example of both cases is given in Figure~\ref{fig:4.1}. From Definition~\ref{Defn 4.1}, we can distinguish four main types of partitions $\Omega$, illustrated in Figure~\ref{fig:4.2}, corresponding to
\begin{enumerate}
    \item Infinite number of Fourier supports, with no rays, which will be denoted by $\Omega^I$: $n_m = -\infty$ and $n_M = \infty$ (i.e. $n \in \Z$),
    \item Infinite number of Fourier supports, with a left ray, which will be denoted by $\Omega_{Lray}^{I}$: $n_m$ finite, $\nu_{n_m} = -\infty$, and $n_M = \infty$, i.e. $\Omega_{n_m} = (-\infty , \nu_{n_m +1}]$,
    \item Infinite number of Fourier supports, with a right ray, which will be denoted by $\Omega_{Rray}^{I}$: $n_m = -\infty$, $n_M$ finite, and $\nu_{n_M} = \infty$, i.e. $\Omega_{n_M -1} = [\nu_{n_M -1} , \infty)$,
    \item Finite number of Fourier supports, with both rays, which will be denoted by $\Omega_{rays}^{F}$: both $n_m$ and $n_M$ are finite, $\nu_{n_m} = -\infty$ and $\nu_{n_M} = \infty$, i.e. $\Omega_{n_m} = (-\infty , \nu_{n_m +1}]$ and $\Omega_{n_M -1} = [\nu_{n_M -1} , \infty)$.\\
\end{enumerate}

We define the length of a Fourier support by $|\Omega_n| = \nu_{n+1} - \nu_n$, with straightforward adaptions for $\Omega^{*}$. 

\begin{figure}[t!]
    \centering
    \includegraphics[width = \linewidth]{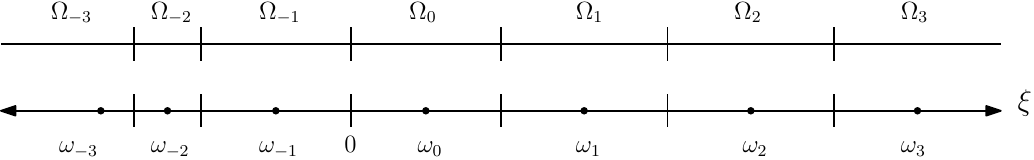}
    \caption{Example of the supports centers in the case of $\Omega^{F}_{rays}$.}
    \label{fig:4.3}
\end{figure}

When $\Omega_n$ is compact, its center is defined by $\omega_n =(\nu_n + \nu_{n+1})/2$, and when rays occur, we symmetrize their adjacent support center, i.e.
\begin{align*}
\omega_{n_m} &= \nu_{n_m +1} - \frac{|\Omega_{n_m +1}|}{2} = \frac{3 \nu_{n_m +1} - \nu_{n_m +2}}{2},  \\
\omega_{n_M -1} &= \nu_{n_M -1} + \frac{|\Omega_{n_M -2}|}{2} = \frac{3 \nu_{n_M -1} - \nu_{n_M -2}}{2}. 
\end{align*}
An example is given in Figure~\ref{fig:4.3} for the case of $\Omega_{rays}^F$. The adaption for $\Omega^{*}$ is straightforward except for $\omega_{-1}$, where
$\omega_{-1} = (\nu_{-1} + \nu_1)/2$.
For the remainder of the paper, we will be referring to sets indexed by $\N$, however, all definitions and results can be restated with sets indexed by $\N^{*}$. We are now ready to define an Empirical Wavelet System.

% ===================================================================
%=================================================================
\subsection{Empirical Wavelet Systems}
%=================================================================
We first define the construction of Empirical Wavelet Systems (EWS). We introduce a new condition, Property 2 in Definition~\ref{def:ews}, which wasn't given in the original work in \cite{cewt}, but will be useful to manage cases when the mother wavelet will not be of compact support. Unless specified, all statements are given for a partition $\Omega$, but can easily be reformulated in the case of $\Omega^*$.\\
\begin{definition}\label{def:ews}
Given a partition $\Omega$ of the frequency domain, let $\psi \in L^2(\R)$ be a function such that its Fourier transform $\widehat{\psi}$ has the following two properties:
\begin{enumerate}
    \item $\widehat{\psi}$ is localized around the zero frequency,   
    \item There exists a subset $E \subseteq \supp{\widehat{\psi}}$ and $0 \leq \delta < 1$, such that
    $$\int_{E} |\widehat{\psi}(\xi)|^2 d\xi = (1-\delta)||\widehat{\psi}||_{L^2}^{2}.$$
    This property guarantees that $\psi$ is mostly supported by $E$.
\end{enumerate}
An \textbf{empirical wavelet system} generated by $\psi$ is the collection
$\{ \psi_{n,b} \mid n \in \N, b \in \R \},$ which will be defined either in the frequency domain by
\begin{align}
\forall \xi \in \R, \; \widehat{\psi}_{n,b}(\xi) = E_{-b} T_{\omega_n} D_{a_n} \widehat{\psi}(\xi) = e^{-2\pi i b\xi} \left( |a_{n}|^{-\frac{1}{2}} \widehat{\psi} \left( \frac{\xi - \omega_n}{a_n} \right) \right),    \label{ewsfreq}
\end{align}
or the time domain by 
\begin{align}
\forall t \in \R, \; \psi_{n,b}(t) = T_b E_{\omega_n} D_{\frac{1}{a_n}} \psi(t) = e^{2\pi i \omega_n (t - b)} \left( |a_{n}|^{\frac{1}{2}} \psi \left( a_n (t - b) \right) \right),     \label{ewstime} 
\end{align}
where $\omega_n$ is the center of the Fourier support $\Omega_n$ and $a_n \in \R \setminus \{0\}$ is a scaling factor whose choice depends on $\Omega_n$ and $\widehat{\psi}$.
\end{definition}
Let us make a few remarks about the definition above. First, observe that $\omega_n$ and $a_n$ are discrete parameters and $b$ is a continuous parameter. Second, the parameters $\omega_n$ and $a_n$ are irregularly sampled but linked to each other.

The first property in the definition implies that each $\widehat{\psi}_n$ will be localized around the center $\omega_n$. The scaling factor $a_n$ has to be chosen such that if $\supp{\widehat{\psi}}_{n}$ is compact, then it should have a width of about $|\Omega_n|$ to ensure $\widehat{\psi}_{n}$ is mostly localized on $\Omega_n$. It must also be chosen in such a way that only consecutive wavelet Fourier supports intersect. Let us denote $\supp{\widehat{\psi}} = S$ and $\supp{\widehat{\psi}_n} = S_n$. We have then $S_n = a_n S + \omega_n$, or equivalently $S = a^{-1}_{n} (S_n - \omega_n)$ and $|S_n| = |a_n| |S|$. From this observation, if $S$ is compact, then we define $a_n$ by
\begin{align}
a_n = \frac{|S_n|}{|S|},  \label{compactsupp an}
\end{align}
so that $\Omega_n \subseteq S_n$ and $|\Omega_n| \leq |S_n|$, for compact $\Omega_n$. When $S$ is not compact, we take $E \subset S$ as described in property $2$ and define $a_n$ by
\begin{align}
a_n = \frac{|\Omega_n|}{|E|},   \label{notcompactsupp an}
\end{align}
for compact $\Omega_n$. In the case of the left and right rays $\Omega_{n_m}$ and $\Omega_{n_M - 1}$, respectively, we define $a_{n_m}$ and $a_{n_M - 1}$ in the following way:\\

\begin{enumerate}
\item If $S = \R$, then $a_{n_m} = \frac{2(\nu_{n_m + 1} - \omega_{n_m})}{|E|}$ and $a_{n_M - 1} = \frac{2(\omega_{n_M - 1} - \nu_{n_M - 1})}{|E|}$,      %\label{suppRRan}   \\
\item If $S = [\xi_L, \infty)$, then $a_{n_m} = -\frac{2(\nu_{n_m + 1} - \omega_{n_m})}{|E|}$ and $a_{n_M - 1} = \frac{2(\omega_{n_M - 1} - \nu_{n_M - 1})}{|E|}$,
%\label{suppRrayan}   \\
\item If $S = (-\infty, \xi_H]$, then $a_{n_m} = \frac{2(\nu_{n_m + 1} - \omega_{n_m})}{|E|}$ and $a_{n_M - 1} = -\frac{2(\omega_{n_M - 1} - \nu_{n_M - 1})}{|E|}$.\\
%\label{suppLrayan}
\end{enumerate}
Noticing that Property 2 in Definition~\ref{def:ews} implies that 
$$\int_{\R \setminus E} |\widehat{\psi}(\xi)|^2 d\xi = \delta ||\widehat{\psi}||^{2}_{L^2},$$
which, combined with Property 1, gives that $\forall \epsilon > 0$, $\exists$ $E \subseteq \R$, about $\xi = 0$, such that 
$$\int_{\R \setminus E} |\widehat{\psi}(\xi)|^2 d\xi < \epsilon.$$
Hence, if $E + r$, with $r >0$, is any enlargement of $E$, then 
$$\lim_{r \rightarrow \infty} \int_{\R \setminus (E +r)} |\widehat{\psi}(\xi)|^2 d\xi = 0$$
It then directly follows that $\widehat{\psi}_n$ inherits the same property. Therefore, an EWS corresponds to a set of bandpass filters associated to the given partition $\Omega$.

For notational convenience, for the remainder of the paper we will denote $\psi_n = E_{\omega_n} D_{\frac{1}{a_n}} \psi$ and $\widehat{\psi}_n = T_{\omega_n} D_{a_n} \widehat{\psi}$. Thus, $\psi_{n,b} = T_b \psi_n$ and $\widehat{\psi}_{n,b} = E_{-b} \widehat{\psi}_n$. \\
Equipped with an EWS, we can now define the continuous empirical wavelet transform (CEWT) by\\
\begin{definition}
Let $f \in L^2(\R)$. The \textbf{continuous empirical wavelet transform} of $f$ is defined by
\begin{align}\label{cewt}
(\mathcal{E}_{\psi}f)(n,b) = \langle \widehat{f}, E_{-b}\widehat{\psi}_n \rangle = \langle f,T_{b}\psi_n \rangle.   
\end{align}
\end{definition}
Similarly to classic wavelets, countable families can be obtained by sampling the translation parameter $b$ from Definition~\ref{def:ews}. This leads to the definition of discrete empirical wavelet system (DEWS).\\

\begin{definition}\label{def:dews}
Let $b=kb_n$, where $k\in\Z$, and $\{ b_n \}_{n \in \N} \subset \R \setminus \{ 0 \}$. Then a discrete empirical wavelet system (DEWS) will correspond to the family of functions, $\forall n\in\N,\forall k\in\Z$, either defined in the Fourier domain by,
$$\forall \xi \in \R, \; \widehat{\psi}_{n,k}(\xi) = E_{-b_{n}k} T_{\omega_n} D_{a_n} \widehat{\psi}(\xi) = e^{-2\pi i b_{n}k\xi} \left( |a_{n}|^{-\frac{1}{2}} \widehat{\psi} \left( \frac{\xi - \omega_n}{a_n} \right) \right),$$
or in the time domain by,
$$\forall t \in \R, \; \psi_{n,k}(t) = T_{b_{n}k} E_{\omega_n} D_{\frac{1}{a_n}} \psi(t) = e^{2\pi i \omega_n (t - b_{n}k)} \left( |a_{n}|^{\frac{1}{2}} \psi \left( a_n (t - b_nk) \right) \right).$$
The corresponding discrete empirical wavelet transform (DEWT) is then given by
$$(\mathcal{E}_{\psi}f)(n,k) = \langle \widehat{f}, E_{-b_{n}k}\widehat{\psi}_n \rangle = \langle f,T_{b_{n}k}\psi_n \rangle.$$
\end{definition}

One final remark concerning an EWS. Notice that if the given partition is either $\Omega^{I}_{Lray}$, $\Omega^{I}_{Rray}$ or $\Omega^{F}_{rays}$, and $\supp{\widehat{\psi}}$ is compact, then there does not exist functions $\widehat{\psi}_{n_m}$ and $\widehat{\psi}_{n_M - 1}$ whose supports can cover $\Omega_{n_m}$ and $\Omega_{n_M - 1}$, respectively. Therefore, with these three partitions, we only consider the subset of $\R$ that only contains the compact Fourier supports. That is, if the given partition is

\begin{enumerate}
    \item  $\Omega^{I}_{Lray}$, then $\Gamma_{Lray} = \displaystyle\bigcup_{n \in \N \setminus \{n_m\}} \Omega_n$,% and the space $L^2(\Gamma_{Rray})$,
    \item $\Omega^{I}_{Rray}$, then $\Gamma_{Rray} = \displaystyle\bigcup_{n \in \N \setminus \{n_M - 1\}} \Omega_n$,% and the space $L^2(\Gamma_{Lray})$,
    \item $\Omega^{F}_{rays}$, then $\Gamma_{C} = \displaystyle\bigcup_{n \in \N \setminus \{n_m, n_{M} - 1\}} \Omega_n$,% and the space $L^2(\Gamma_{C})$.
\end{enumerate}
and therefore, we introduce the spaces $L^2_{Lray}(\R), L^2_{Rray}(\R)$ and $L^2_C(\R)$.\\

\begin{definition}
 The spaces $L^2_{Lray}(\R), L^2_{Rray}(\R)$ and $L^2_C(\R)$ are the subspaces of $L^2(\R)$ defined by
 $$L^2_X(\R)=\{f\in L^2(\R) \rvert\; \supp \hat{f}\subseteq \Gamma_X\}\qquad \text{where}\quad X\in\{Lray,Rray,C\}.$$
\end{definition}

%% ===================================================================
\section{Frames of wave packet systems}\label{sec3}

In this section, we recall the notion of wave packet systems and some of their properties which will be useful to prove our results in the next section. Wave packet systems are defined by\\

\begin{definition}
  Let $\Psi=\{\psi^l: 1\leq l\leq L\}\subset L^2(\R^n)$, where $L$ is a finite integer, and let $S\subseteq GL_n(\R)\times\R^n$. The continuous wave packet system relative to $S$ generated by $\Psi$ is the collection defined by,
  \begin{equation}
    \W_S^{(i)}(\Psi)=\{U^{(i)}_{(a,\nu,y)}\psi^l:(a,\nu)\in S,y\in\R^n,1\leq l\leq L\},
  \end{equation}
  where the operators $U^{(i)},1\leq i\leq 5$ are defined by
  \begin{align*}
    U^{(1)}_{(a,\nu,y)}&= D_aT_yE_\nu, & U^{(2)}_{(a,\nu,y)}&= T_yD_aE_\nu, \\
    U^{(3)}_{(a,\nu,y)}&= E_\nu D_a T_y, & U^{(4)}_{(a,\nu,y)}&= T_yE_\nu D_a, \\
    U^{(5)}_{(a,\nu,y)}&= E_\nu T_yD_a. & 
  \end{align*}
\end{definition}

In \cite{Labate2004}, the authors define the notion of continuous Parseval frame wave packet for the important case of wave packet systems associated with a reproducing formula.\\

\begin{definition}\label{def:contparsevalframe}
  Let $S\subseteq GL_n(\R)\times \R^n$, and $\lambda$ be a measure on $S$. Let $1\leq i\leq 5$, the system $\W_S^{(i)}(\Psi)$ is a continuous Parseval frame wave packet system relative to $(S,\lambda)$ for $L^2(\R^n)$, provided that the functions $(a,\nu,y)\mapsto\langle f,U^{(i)}_{(a,\nu,y)}\psi^l\rangle$ are $\lambda-$measurable for all $f,\psi\in L^2(\R^n)$ and 
  $$\|f\|_{L^2}^2=\sum_{l=1}^L\int_{S\times\R^n} |\langle f,U^{(i)}_{(a,\nu,y)}\psi^l\rangle|^2 d\lambda(a,\nu)dy,$$
  for all $f\in L^2(\R^n)$.\\
\end{definition}

In most applications, frames are built from countable families of functions. A general formalism for countable wave packet systems have been introduced in \cite{unifiedii} (see also \cite{HERNANDEZ2004111}, and \cite{Labate2004}). We, hereafter, recall such definition as well as two important results which give us some conditions for when countable wave packet systems are frames, as well as to help characterize their frame bounds. Even though our work will only consider the one-dimensional case, we keep stating these results in the $n-$dimensional case as they are given in \cite{unifiedii}.

Let $\mathcal{P}$ be a countable indexing set, $\{ g_p \}_{p \in \mathcal{P}}$ be a collection of functions in $L^2(\R^n)$ and $\{ C_p \}_{p \in \mathcal{P}}$ be an associated collection of matrices in $GL_n(\R)$ corresponding to deformations of the functions $g_p$. A corresponding sampling of $\R^n$ is given by

$$\Lambda = \bigcup_{p \in \mathcal{P}} C^{I}_{p} \Z^n,$$
where $C^{I}_{p} = (C^{T}_{p})^{-1}$. Given $\alpha \in \Lambda$, we define
$$\mathcal{P}_{\alpha} = \{ p \in \mathcal{P} \mid \alpha \in C^{I}_{p} \Z^n \} = \{ p \in \mathcal{P} \mid C^{T}_{p} \alpha \in \Z^n \}.$$
Notice that $\mathcal{P}_{\alpha} \subset \mathcal{P}$, and that if $\alpha = 0 \in \Lambda$, then $\mathcal{P}_0 = \mathcal{P}$. The following set
\begin{equation}
\mathcal{D} = \{ f \in L^2(\R^n) \mid \hat{f} \in L^{\infty}(\R^n) \text{ and } \supp{\hat{f}} \text{ is compact} \}, \label{dense set}
\end{equation}
is dense in $L^2(\R^n)$. Moreover, for any Lebesque measurable subset $X \subseteq \R$, if we define

$$\Check{L}^2(X) = \{ f \in L^2(\R) \mid \supp{\hat{f}} \subset X \},$$ %\label{lineaer subspace} \\
then
\begin{equation}
\mathcal{D}' = \{ f \in L^2(\R) \mid \hat{f} \in L^{\infty}(\R) \text{ and } \supp{\hat{f}} \subset X \text{ is compact} \}, \label{dense set X}
\end{equation}
is dense in $\check{L}^2(X)$, \cite{christensen2016explicit}. The following theorem \cite{unifiedii} shows that generalized shift-invariant systems, i.e. families of the form
\begin{equation}
\{ T_{C_{p}k} g_{p} \mid p \in \mathcal{P}, k \in \Z^n \}, \label{oversampled} 
\end{equation}
are Parseval frames for $L^2(\R^n)$.\\
\begin{theorem} \label{PF thm}
Let $\mathcal{P}$ be a countable indexing set, $\{ g_p \}_{p \in \mathcal{P}}$ a collection of functions in $L^2(\R^n)$ and $\{ C_p \}_{p \in \mathcal{P}} \subset GL_n(\R)$. Suppose that 
\begin{align}
L(f) = \sum_{p \in \mathcal{P}} \sum_{m \in \Z^n} \int_{\supp{\hat{f}}} |\hat{f}(\xi + C^{I}_{p}m)|^2 \frac{1}{|\det C_p|} |\hat{g}_{p}(\xi)|^2 d\xi < \infty, \label{LIC}
\end{align}
for all $f \in \mathcal{D}$, where $C^{I}_{p} = (C^{T}_{p})^{-1}$. Then systems of the form of \eqref{oversampled} are Parseval frames for $L^2(\R^n)$ if and only if 
$$\sum_{p \in \mathcal{P}_{\alpha}} \frac{1}{|\det C_p|} \hat{g}_{p}(\xi) \overline{\hat{g}_{p}(\xi + \alpha)} = \delta_{\alpha,0} \;\; for \; a.e. \; \xi \in \R^n,$$
for each $\alpha \in \Lambda$, where $\delta_{\alpha,0}$ is the Kronecker Delta function for $\R^n$.\\
\end{theorem}

We will also need the following theorem established in \cite{Labate2004} and improved upon in \cite{Christensen2008a}.

\begin{theorem} \label{frame thm}
Let $\{g_p\}_{p \in \mathcal{P}} \subset L^2(\R^n)$, $\mathcal{P}$ be a countable index set, and $\{C_p\}_{p \in \mathcal{P}} \subset GL_n(\R)$. If 
$$B = \sup_{\xi \in \R^n} \sum_{p \in \mathcal{P}}\sum_{k \in \Z^n} \frac{1}{|\det C_p|} |\hat{g}_{p}(\xi)\hat{g}_{p}(\xi - C^{I}_{p}k)| < \infty,$$
then the system $\{ T_{C_{p}k} g_{p} \mid p \in \mathcal{P}, k \in \Z^n \}$ is a Bessel sequence with bound $B$. Furthermore, if 
$$A = \inf_{\xi \in \R^n} \bigg( \sum_{p \in \mathcal{P}}\frac{1}{|\det C_p|}|\hat{g}_{p}(\xi)|^2 - \sum_{p \in \mathcal{P}}\sum_{k \neq 0}\frac{1}{|\det C_p|} |\hat{g}_{p}(\xi)\hat{g}_{p}(\xi - C^{I}_{p}k)| \bigg) > 0,$$
then $\{ T_{C_{p}k} g_{p} \mid p \in \mathcal{P}, k \in \Z^n \}$ is a frame for $L^2(\R_n)$ with bounds $A$ and $B$.
\end{theorem}

%% =====================================================================
\section{Frames of empirical wavelets}\label{sec4}
In this section, we investigate the possibility of building EW frames. We first provide a condition on the EW family to be a continuous Parseval wave packet frame. Then, we study the construction of discrete frames, addressing two cases: when the mother wavelet Fourier transform has a compact support or not.

\subsection{Empirical wavelets as a continuous Parseval wave packet frame}
The next theorem provides a condition to guarantee that an EW family $\{\psi_{b,n}\}$, as defined in Definition~\ref{def:ews}, forms a Parseval wave packet frame (see Definition~\ref{def:contparsevalframe}).\\

\begin{theorem}
  The empirical wavelet system $\{\psi_{b,n}\}$ is a continuous Parseval frame wave packet system relative to $(\R,\Z)$ (or $(\R,\Z^*)$) for $L^2(\R)$ if and only if
  $$\sum_{n\in \Z}\frac{1}{a_n}\left|\hat{\psi}\left(\frac{\xi-\omega_n}{a_n}\right)\right|^2=1\qquad\text{for a.e. }\; \xi\in\R.$$
\end{theorem}

\proof
Since for all $f\in L^2(\R),|\langle f,T_bE_{\omega_n}D_{1/a_n}\psi\rangle|=|\langle \hat{f},E_{-b}T_{\omega_n}D_{a_n}\hat{\psi}\rangle|$, using Plancherel theorem, we have
\begin{align*}
  \int_\R |\langle f,T_bE_{\omega_n}D_{1/a_n}\psi\rangle|^2db&= \int_\R |\langle \hat{f},E_{-b}T_{\omega_n}D_{a_n}\hat{\psi}\rangle|^2db\\
  &=\int_\R \left|\int_\R \hat{f}(\xi) \overline{E_{-b}T_{\omega_n}D_{a_n}\hat{\psi}(\xi)}d\xi\right|^2 db\\
  &=\int_\R \left|\int_\R \hat{f}(\xi) \frac{1}{\sqrt{a_n}} \overline{\hat{\psi}\left(\frac{\xi-\omega_n}{a_n}\right)} e^{2\imath\pi b\xi} d\xi\right|^2 db\\
  &=\frac{1}{a_n}\int_\R \left|\left(\hat{f}\cdot\overline{\hat{\psi}\left(\frac{\cdot-\omega_n}{a_n}\right)}\right)^\vee (b)\right|^2 db\\
  &=\frac{1}{a_n}\int_\R |\hat{f}(\xi)|^2\left|\hat{\psi}\left(\frac{\xi-\omega_n}{a_n}\right)\right|^2 d\xi.
\end{align*}
Now, if we sum with respect $n$, we get
\begin{align*}
  \sum_{n\in\Z}\int_\R |\langle f,T_bE_{\omega_n}D_{1/a_n}\psi\rangle|^2db&=\sum_{n\in\Z}\frac{1}{a_n}\int_\R |\hat{f}(\xi)|^2\left|\hat{\psi}\left(\frac{\xi-\omega_n}{a_n}\right)\right|^2 d\xi\\
  &=\int_\R |\hat{f}(\xi)|^2 \sum_{n\in\Z}\frac{1}{a_n}\left|\hat{\psi}\left(\frac{\xi-\omega_n}{a_n}\right)\right|^2 d\xi,
\end{align*}
therefore, if $\sum_{n\in\Z}\frac{1}{a_n}\left|\hat{\psi}\left(\frac{\xi-\omega_n}{a_n}\right)\right|^2=1$ for a.e. $\xi\in\R$, we get 
$$\sum_{n\in\Z}\int_\R |\langle f,T_bE_{\omega_n}D_{1/a_n}\psi\rangle|^2db =\int_\R |\hat{f}(\xi)|^2d\xi=\|\hat{f}\|^2=\|f\|^2.$$
Thus, we obtain that the empirical wavelet system $\{\psi_{b,n}\}$ fulfills the conditions of Definition~\ref{def:contparsevalframe}, and is, therefore, a continuous Parseval frame wave packet.
\endproof

\subsection{Empirical wavelets as discrete wave packet frames}
In this section, we consider countable families of EW as defined in Definition~\ref{def:dews}. We distinguish two cases: when the Fourier transform of the mother wavelet $\psi$ has compact support or not. For each case, we provide the conditions to have the corresponding discrete empirical wavelet system (DEWS) to be a discrete wave packet frame.

%=================================================================
\subsubsection{When $\widehat{\psi}$ has compact support}
%=================================================================
In this section, we prove Theorem~\ref{result1} and \ref{th:dewswithrays} which provide a necessary and sufficient condition for a DEWS associated to a partition $\Omega$ to be a Parseval frame. We first investigate the case of partitions without rays.

\begin{theorem} \label{result1}
Let $\Omega^{I}$ be a given partition, let $\psi \in L^2(\R)$ be such that $\widehat{\psi}$ has compact support, and $\N$ be a countable indexing set. Then the system $\{ T_{b_{n}k}\psi_{n} \mid n \in \N, k \in \Z  \}$ is a Parseval frame for $L^2(\R)$ if and only if
$$\sum_{n \in \mathcal{N_\alpha}} \frac{1}{|b_n||a_n|} \widehat{\psi} \left( \frac{\xi - \omega_n}{a_n} \right) \overline{\widehat{\psi} \left( \frac{\xi + \alpha - \omega_n}{a_n} \right)} = \delta_{\alpha, 0} \quad\;for \; a.e. \quad \xi \in \R,$$
for each $\alpha \in \Lambda$, where

\begin{gather}
\Lambda = \bigcup_{n \in \N} b^{-1}_{n}\Z,  \\
\N_{\alpha} = \{ n \in \N \mid \alpha \in b^{-1}_{n} \Z \} = \{ n \in \N \mid b_{n} \alpha \in \Z \},
\end{gather}
and $\delta_{\alpha,0}$ is the Kronecker Delta function on $\R$, where $\delta_{\alpha, 0} = 1$, if $\alpha = 0$ and $\delta_{\alpha, 0} = 0$, if $\alpha \neq 0$.
\end{theorem}

\proof
We aim at applying Theorem~\ref{PF thm} in the 1-dimensional case, to do so, we need to show that \eqref{LIC} holds. Letting $\mathcal{P} = \N$, $g_p = \psi_n = E_{\omega_n} D_{\frac{1}{a_n}} \psi$ and $C_p = b_n$, based onto \eqref{LIC} we have, 

\begin{align*}
  L(f) &= \sum_{n \in \N}\sum_{m \in \Z} \int_{\supp{\hat{f}}}\left|\hat{f}(\xi + b^{-1}_{n} m)\right|^2 \frac{1}{|b_n||a_n|} \left|\widehat{\psi} \left( \frac{\xi - \omega_n}{a_n} \right)\right|^2 d\xi \\
  &= \sum_{n \in \N}\sum_{m \in \Z} \int_{\supp{\hat{f}} \cap \supp{\widehat{\psi}_n}} \left|\hat{f}(\xi + b^{-1}_{n} m)\right|^2 \frac{1}{|b_n||a_n|} \left|\widehat{\psi} \left( \frac{\xi - \omega_n}{a_n} \right)\right|^2 d\xi.
\end{align*}
Since $\supp{\hat{f}}$ is compact and $\supp{\widehat{\psi}_n}$ is compact, for each fixed $n \in \N$, there are finitely many $m \in \Z$ (say, $M$ of them) such that $L(f)$ is nonzero. And since $||\hat{f}||_{\infty} < \infty$, we have that 

$$L(f) \leq M ||\hat{f}||^{2}_{\infty} \sum_{n \in \N} \int_{\supp{\hat{f}} \cap \supp{\widehat{\psi}_n}} \frac{1}{|b_n||a_n|} \left| \widehat{\psi}\left( \frac{\xi - \omega_n}{a_n} \right)  \right|^2 d\xi.$$
Furthermore, there are finitely many $\supp{\widehat{\psi}_n}$ intersecting $\supp{\hat{f}}$. Thus,
$$\sum_{n \in \N} \int_{\supp{\hat{f}} \cap \supp{\widehat{\psi}_n}} \frac{1}{|b_n||a_n|} \left| \widehat{\psi}\left( \frac{\xi - \omega_n}{a_n} \right)  \right|^2 d\xi<\infty,$$
which implies that $L(f) < \infty$. This proves that \eqref{LIC} holds. The result now follows from applying Theorem~\ref{PF thm} in the 1-dimensional case.
\endproof

Next, we address the case when the partition has some rays. Note, from Theorem~\ref{result1}, for an EWS, the condition \eqref{LIC} is equivalent to

\begin{align}
L(f) &= \sum_{n \in \N}\sum_{m \in \Z} \int_{\supp{\hat{f}}} \left|\hat{f}(\xi + b^{-1}_{n} m)\right|^2 \frac{1}{|b_n||a_n|} \left|\widehat{\psi} \left( \frac{\xi - \omega_n}{a_n} \right)\right|^2 d\xi \notag\\
&= \sum_{n \in \N}\sum_{m \in \Z} \int_{\supp{\hat{f}} \cap \supp{\widehat{\psi}_n}} \left|\hat{f}(\xi + b^{-1}_{n} m)\right|^2 \frac{1}{|b_n||a_n|} \left|\widehat{\psi} \left( \frac{\xi - \omega_n}{a_n} \right)\right|^2 d\xi < \infty. \label{LICews}
\end{align}

Therefore, for the remaining three types of partitions (i.e. with rays), it will suffice to show that \eqref{LICews} holds. We address the cases of rays, either a single left or right ray, or both left and right rays, by proving the three following theorems. These results show that EW frames can be obtained for the subspaces $L_{Lray}^2(\R), L_{Rray}^2(\R)$ or $L_F^2(\R)$, respectively.

\begin{theorem}\label{th:dewswithrays}
  Let $\Omega^{I}_{Lray}$ be the given partition, $\psi \in L^2(\R)$ be such that $\widehat{\psi}$ has compact support, and $\N$ be a countable indexing set. Then the system $\{ T_{b_{n}k}\psi_{n} \mid n \in \N \setminus \{n_m\}, k \in \Z  \}$ is a Parseval frame for $L_{Lray}^2(\R)$ if and only if
$$\sum_{n \in \N_\alpha} \frac{1}{|b_n||a_n|} \widehat{\psi} \left( \frac{\xi - \omega_n}{a_n} \right) \overline{\widehat{\psi} \left( \frac{\xi + \alpha - \omega_n}{a_n} \right)} = \delta_{\alpha, 0} \quad for \; a.e. \quad \xi \in \Gamma_{Lray},$$
for each $\alpha \in \Lambda$, where
\begin{gather}
\Lambda = \bigcup_{n \in \N \setminus \{n_m\}} b^{-1}_{n}\Z,  \\
\N_\alpha = \{ n \in \N \setminus \{n_m\} \mid \alpha \in b^{-1}_{n} \Z \} = \{ n \in \N \setminus \{n_m\} \mid b_{n} \alpha \in \Z \},
\end{gather}
and $\delta_{\alpha,0}$ is the Kronecker Delta function on $\R$.
\end{theorem}
\proof
We show that \eqref{LICews} holds for $f \in \mathcal{D}'$ as defined in \eqref{dense set X}, where $X = \Gamma_{Lray}$. We argue as above: Since $\supp{\hat{f}}$ and $\supp{\widehat{\psi}_n}$ are both compact, for each fixed $n \in \N \setminus \{n_m\}$, there are finitely many $m \in \Z$ such that $L(f)$ is nonzero. And since finitely many $\supp{\widehat{\psi}_n}$ intersect $\supp{\hat{f}}$, together with the fact that $\|\hat{f}\|_{\infty} < \infty$, we see that $L(f) < \infty$. The result follows from Theorem~\ref{PF thm}.
\endproof

\begin{theorem}
  Let $\Omega^{I}_{Rray}$ be the given partition, $\psi \in L^2(\R)$ be such that $\widehat{\psi}$ has compact support, and $\N$ be a countable indexing set. Then the system $\{ T_{b_{n}k}\psi_{n} \mid n \in \N \setminus \{n_M - 1\}, k \in \Z  \}$ is a Parseval frame for $L_{Rray}^2(\R)$ if and only if
$$\sum_{n \in \mathcal{N_\alpha}} \frac{1}{|b_n||a_n|} \widehat{\psi} \left( \frac{\xi - \omega_n}{a_n} \right) \overline{\widehat{\psi} \left( \frac{\xi + \alpha - \omega_n}{a_n} \right)} = \delta_{\alpha, 0} \quad for \; a.e. \quad \xi \in \Gamma_{Rray},$$
for each $\alpha \in \Lambda$, where

\begin{align}
\Lambda &= \bigcup_{n \in \N \setminus \{n_M -1\}} b^{-1}_{n}\Z,  \\
\N_{\alpha} &= \{ n \in \N \setminus \{n_M -1\} \mid \alpha \in b^{-1}_{n} \Z \} = \{ n \in \N \setminus \{n_M -1\} \mid b_{n} \alpha \in \Z \},
\end{align}
and $\delta_{\alpha,0}$ is the Kronecker Delta function on $\R$.
\end{theorem}

\proof
Take $\mathcal{D}'$ as defined in \eqref{dense set X}, where $X = \Gamma_{Rray}$. Then the result follows from a similar argument as used in the proof of Theorem~\ref{th:dewswithrays}.
\endproof

\begin{theorem}
  Let $\Omega^{F}_{rays}$ be the given partition, $\psi \in L^2(\R)$ be such that $\widehat{\psi}$ has compact support, and $\N$ be a countable indexing set. Then the system $\{ T_{b_{n}k}\psi_{n} \mid n \in \N \setminus \{n_m,n_M - 1\}, k \in \Z  \}$ is a Parseval frame for $L_F^2(\R)$ if and only if
$$\sum_{n \in \mathcal{N_\alpha}} \frac{1}{|b_n||a_n|} \widehat{\psi} \left( \frac{\xi - \omega_n}{a_n} \right) \overline{\widehat{\psi} \left( \frac{\xi + \alpha - \omega_n}{a_n} \right)} = \delta_{\alpha, 0} \quad for \; a.e. \quad \xi \in \Gamma_{C},$$
for each $\alpha \in \Lambda$, where
\begin{gather}
\Lambda = \bigcup_{n \in \N \setminus \{n_m, n_M -1\}} b^{-1}_{n}\Z,  \\
\N_{\alpha} = \{ n \in \N \setminus \{n_m, n_M -1\} \mid \alpha \in b^{-1}_{n} \Z \} = \{ n \in \N \setminus \{n_m, n_M -1\} \mid b_{n} \alpha \in \Z \},
\end{gather}
and $\delta_{\alpha,0}$ is the Kronecker Delta function on $\R$.
\end{theorem}
\proof
Take $\mathcal{D}'$ as defined in \eqref{dense set X}, where $X = \Gamma_{C}$. Then the result follows from a similar argument as used in the proof of Theorem~\ref{th:dewswithrays}.
\endproof

%=================================================================
\subsubsection{When $\widehat{\psi}$ has non-compact support}
%=================================================================
In this section, we aim to prove equivalent results like Theorem~\ref{frame thm} to get sufficient conditions for an EWS associated to any of the four main partitions $\Omega$ to be a frame, in the case when $\hat{\psi}$ does not have a compact support.

\begin{theorem}
Let $\Omega$ be any of the four main partitions, $\psi \in L^2(\R)$ be such that $\widehat{\psi}$ has non-compact support, and $\N$ be a countable indexing set.
If
$$B = \sup_{\xi \in \R} \sum_{n \in \N}\sum_{k \in \Z} \frac{1}{|b_n| |a_n|} \left| \widehat{\psi} \left( \frac{\xi - \omega_{n}}{a_n} \right) \widehat{\psi} \left( \frac{\xi - b^{-1}_{n}k - \omega_n}{a_n} \right) \right| < \infty,$$ 
and 
\begin{align*}
A = \inf_{\xi \in \R} &\left( \sum_{n \in \N} \frac{1}{|b_n| |a_n|}  
\left|  \widehat{\psi} \left( \frac{\xi - \omega_n}{a_n} \right) \right|^2 \right. \\
&-\left. \sum_{n \in \N}\sum_{k \neq 0} \frac{1}{|b_n| |a_n|} \left| \widehat{\psi} \left( \frac{\xi - \omega_n}{a_n} \right) \widehat{\psi} \left( \frac{\xi - b^{-1}_{n}k - \omega_n}{a_n} \right) \right| \right) > 0,
\end{align*}
then the system $\{ T_{b_{n}k}\psi_n \mid n \in \N, k \in \Z  \}$, associated with the partition $\Omega$, is a frame for $L^2(\R)$ with frame bounds $A$ and $B$.
\end{theorem}

\proof
It suffices to prove the theorem for $f \in \mathcal{D}$, as defined in \eqref{dense set}. Let $n \in \N$, by Plancherel's theorem, we have
$$\sum_{k \in \Z}|\langle f,T_{b_{n}k}\psi_n\rangle |^2 = \sum_{k \in \Z} \left| \underbrace{\int_{\R} \hat{f}(\xi) \overline{\widehat{\psi}_n(\xi)} e^{2\pi i b_n k \xi} d\xi}_I  \right|^2.$$
Let $\mathbb{T} = [0,1)$, since $\R$ can be written as a disjoint union, i.e. $\R = \bigcup_{l \in \Z} b^{-1}_{n}(\mathbb{T} - l)$, the integral, $I$, above can be written in the form 
$$I=\int_{b^{-1}_{n}\mathbb{T}} \left( \sum_{l \in \Z} \hat{f}(\xi - b^{-1}_{n}l) \overline{\widehat{\psi}_n(\xi - b^{-1}_{n}l)}  \right) e^{2\pi i b_{n}k \xi} d\xi.$$
Thus, we can then write
\begin{align*}
\sum_{k \in \Z}|\langle f,T_{b_{n}k} &\psi_n \rangle |^2 = \sum_{k \in \Z} \left| \int_{\R} \hat{f}(\xi) \overline{\widehat{\psi}_n(\xi)} e^{2\pi i b_n k \xi} d\xi  \right|^2   \\
&= \sum_{k \in \Z} \left| \int_{b^{-1}_{n}\mathbb{T}} \left( \sum_{l \in \Z} \hat{f}(\xi - b^{-1}_{n}l) \overline{\widehat{\psi}_n(\xi - b^{-1}_{n}l)}  \right) e^{2\pi i b_{n}k \xi} d\xi  \right|^2   \\
&= |b_n|^{-2}\sum_{k \in \Z} \left| \frac{1}{|b_n|^{-1}}\int_{b^{-1}_{n}\mathbb{T}} \left( \sum_{l \in \Z} \hat{f}(\xi - b^{-1}_{n}l) \overline{\widehat{\psi}_n(\xi - b^{-1}_{n}l)}  \right) e^{2\pi i \frac{k}{b^{-1}_{n}} \xi} d\xi  \right|^2.
\end{align*}
Since $\sum_{l \in \Z} \hat{f}(\xi - b^{-1}_{n}l) \overline{\widehat{\psi}_n(\xi - b^{-1}_{n}l)}$ is $b^{-1}_{n}\Z$-periodic and belongs to $L^2(b^{-1}_{n}\Z)$ (recall $f \in \mathcal{D}$), by Parseval's identity we obtain
\begin{align*}
\sum_{k \in \Z}|\langle f,T_{b_{n}k}\psi_n\rangle|^2 &= \frac{|b_n|^{-2}}{|b_n|^{-1}} \int_{b^{-1}_{n}\mathbb{T}} \left| \sum_{l \in \Z} \hat{f}(\xi - b^{-1}_{n}l) \overline{\widehat{\psi}_n(\xi - b^{-1}_{n}l)} \right|^2 d\xi    \\
&= \frac{1}{|b_n|} \int_{b^{-1}_{n}\mathbb{T}} \left| \sum_{l \in \Z} \hat{f}(\xi - b^{-1}_{n}l) \overline{\widehat{\psi}_n(\xi - b^{-1}_{n}l)} \right|^2 d\xi   \\
= \frac{1}{|b_n|}& \int_{b^{-1}_{n}\mathbb{T}} \sum_{l \in \Z} \hat{f}(\xi - b^{-1}_{n}l) \overline{\widehat{\psi}_n(\xi - b^{-1}_{n}l)}\; \overline{\sum_{l \in \Z} \hat{f}(\xi - b^{-1}_{n}l) \overline{\widehat{\psi}_n(\xi - b^{-1}_{n}l)}} d\xi \\
= \frac{1}{|b_n|}& \int_{b^{-1}_{n}\mathbb{T}} \sum_{l \in \Z} \hat{f}(\xi - b^{-1}_{n}l) \overline{\widehat{\psi}_n(\xi - b^{-1}_{n}l)}\; \sum_{l \in \Z} \overline{\hat{f}(\xi - b^{-1}_{n}l)} \widehat{\psi}_n(\xi - b^{-1}_{n}l) d\xi.
\end{align*}
By making the change of indices $u=l+k$ and applying the same techniques from \cite{Labate2004}, we obtain
\begin{align*}
&\frac{1}{|b_n|} \int_{b^{-1}_{n}\mathbb{T}} \sum_{l,u \in \Z} \hat{f}(\xi - b^{-1}_{n}l) \overline{\widehat{\psi}_n(\xi - b^{-1}_{n}l)}\;  \overline{\hat{f}(\xi - b^{-1}_{n}u)} \widehat{\psi}_n(\xi - b^{-1}_{n}u) d\xi \\
&= \frac{1}{|b_n|} \int_{b^{-1}_{n}\mathbb{T}} \sum_{l,k \in \Z} \hat{f}(\xi - b^{-1}_{n}l) \overline{\widehat{\psi}_n(\xi - b^{-1}_{n}l)}\;  \overline{\hat{f}(\xi - b^{-1}_{n}(l+k))} \widehat{\psi}_n(\xi - b^{-1}_{n}(l+k)) d\xi  \\
&= \frac{1}{|b_n|} \sum_{k\in\mathbb{Z}} \int_{\mathbb{R}} \hat{f}(\xi) \overline{\widehat{\psi}_n(\xi)} \overline{\hat{f}(\xi - b^{-1}_{n}k)} \widehat{\psi}_n(\xi - b^{-1}_{n}k)   d\xi.
\end{align*}
So that we obtain the following analogous result
$$\sum_{k \in \Z}|\langle f,T_{b_{n}k}\psi_n\rangle|^2 = \frac{1}{|b_n|} \sum_{k \in \Z} \int_{\R} \hat{f}(\xi) \overline{\hat{f}(\xi - b^{-1}_{n}k)}\, \overline{\widehat{\psi}_n(\xi)} \widehat{\psi}_n(\xi - b^{-1}_{n}k) d\xi.  $$
Hence, splitting the terms when $k=0$ and $k\neq 0$ apart from each other, it follows that
\begin{align*}
\sum_{n \in \N}\sum_{k \in \Z} |\langle f, T_{b_{n}k}\psi_n \rangle|^2 = \int_{\R} |\hat{f}(\xi)|^2 \sum_{n \in \N} \frac{1}{|b_n|} |\widehat{\psi}_n(\xi)|^2 d\xi + R(f),
\end{align*}
where
\begin{align*}
R(f) = \sum_{n \in \N}\sum_{k \neq 0} \frac{1}{|b_n|} \int_{\R} \hat{f}(\xi)\overline{\hat{f}(\xi - b^{-1}_{n}k)}\, \overline{\widehat{\psi}_n(\xi)}\widehat{\psi}_n(\xi - b^{-1}_{n}k) d\xi.
\end{align*}
It is also proven in \cite{Labate2004} that
\begin{align*}
R(f) \leq \sum_{n \in \N}\sum_{k \neq 0} \frac{1}{|b_n|} \int_{\R} |\hat{f}(\xi)|^2 |\widehat{\psi}_n(\xi)\widehat{\psi}_n(\xi - b^{-1}_{n}k)| d\xi.
\end{align*}
Furthermore, applying a similar argument as in the proof of Theorem~\ref{frame thm} (see \cite{Christensen2008a}), we get
\begin{align*}
\sum_{n \in \N}\sum_{k \in \Z} |\langle f, & T_{b_{n}k}\psi_n \rangle|^2  \leq \\ &\int_{\R} |\hat{f}(\xi)|^2 \sum_{n \in \N}\sum_{k \in \Z} \frac{1}{|b_n| |a_n|} \left| \widehat{\psi} \left( \frac{\xi - \omega_{n}}{a_n} \right) \widehat{\psi} \left( \frac{\xi - b^{-1}_{n}k - \omega_n}{a_n} \right) \right| d\xi,
\end{align*}
and
\begin{align*}
\sum_{n \in \N}\sum_{k \in \Z} |\langle f, T_{b_{n}k}\psi_n \rangle|^2 \geq &\int_{\R} |\hat{f}(\xi)|^2 \left( \sum_{n \in \N} \frac{1}{|b_n| |a_n|} \left| \widehat{\psi} \left( \frac{\xi - \omega_n}{a_n} \right) \right|^2 \right.\\
&- \left.\sum_{n \in \N}\sum_{k \neq 0} \frac{1}{|b_n| |a_n|} \left| \widehat{\psi} \left( \frac{\xi - \omega_{n}}{a_n} \right) \widehat{\psi} \left( \frac{\xi - b^{-1}_{n}k - \omega_n}{a_n} \right) \right| \right)d\xi.
\end{align*}
By respectively taking the supremum and infimum, we obtain the desired result.
\endproof

\section{Conclusion}\label{sec5}
In this work, we have investigated the possibility to build frames of empirical wavelets. We have provided conditions for the existence of continuous and discrete empirical wavelet frames. We also provided an estimation of the frame bounds in the case of frames based on families of empirical wavelet built from a mother wavelet of non-compact support in the Fourier domain. The provided results will also allow better construction of empirical wavelets for certain class of applications like signal transmission, signal processing, just to mention a few.\\
In terms of future work, similar investigations in higher dimensions would be of great interest as well. However, the presence of some arbitrary geometric constraints in the construction of such empirical wavelets makes the generalization more challenging since no general formalism has been developed yet.

\section*{Acknowledgement}
This work has been sponsored by the Air Force Office of Scientific Research, grant FA9550-21-1-0275.

\end{document}